\documentclass[MSNbibl,nameyear,citesort,seceqn,dvips]{arxbj}

\aid{0}
\volume{19}
\issue{4}
\pubyear{2013}
\firstpage{1404}
\lastpage{1418}
\doi{10.3150/12-BEJSP03} %kopijuoti is PTS

\makeatletter
\newcommand{\dotsim}{\mathop{\sim}^{.}}
\newcommand{\goesd}{\mathop{\longrightarrow}^{\mathcal{L}}}
\makeatother

\begin{document}
\begin{frontmatter}

\title{Aspects of likelihood inference}
\runtitle{Likelihood}

\begin{aug}
%%%% inicialai - be tarpu
\author{\fnms{Nancy} \snm{Reid}\corref{}\ead[label=e1]{reid@utstat.utoronto.ca}\ead[label=u1,url]{www.utstat.utoronto.ca/reid/}}% \and
\runauthor{N. Reid} %% auto
\address{Department of Statistics, University of Toronto,
100 St. George St.,
Toronto, Canada M5S 3G3.\\ \printead{e1,u1}}
\end{aug}

% HISTORY:

% ABSTRACT
%
\begin{abstract}
I review the classical theory of likelihood based inference and
consider how it is being extended and developed for use in complex
models and sampling schemes.
\end{abstract}

% KEYWORDS
% visi is mazosios raides ir pagal abecele
%

\end{frontmatter}

%s1 ###
\section{Introduction}\label{sec1}

Jakob Bernoulli's \emph{Ars Conjectandi} established the field of
probability theory, and founded a long and remarkable mathematical
development of deducing patterns to be observed in sequences of random
events. The theory of statistical inference works in the opposite
direction, attempting to solve the inverse problem of deducing
plausible models from a given set of observations. Laplace pioneered
the study of this inverse problem, and indeed he referred to his method
as that of inverse probability.

The likelihood function, introduced by \citet{fisher22}, puts this
inversion front and centre, by writing the probability model as a
function of unknown parameters in the model. This simple, almost
trivial, change in point of view has profoundly influenced the
development of statistical theory and methods. In the early days,
computing data summaries based on the likelihood function could be
computationally difficult, and various \emph{ad hoc} simplifications
were proposed and studied. By the late 1970s, however, the widespread
availability of computing enabled a parallel development of widespread
implementation of likelihood-based inference. The development of
simulation and approximation methods that followed meant that both
Bayesian and non-Bayesian inferences based on the likelihood function
could be readily obtained.

As a result, construction of the likelihood function, and various
summaries derived from it, is now a nearly ubiquitous starting point
for a great many application areas. This has a unifying effect on the
field of applied statistics, by providing a widely accepted standard as
a starting point for inference.
%(Aspects of the next steps in the inference process continue to be
%widely debated.)

With the explosion of data collection in recent decades, realistic
probability models have continued to grow in complexity, and the
calculation of the likelihood function can again be computationally
very difficult. Several lines of research in active development concern
methods to compute approximations to the likelihood function, or
inference functions with some of the properties of likelihood
functions, in these very complex settings.

In the following section, I will summarize the standard methods for
inference based on the likelihood function, to establish notation, and
then in Section~\ref{sec3} describe some aspects of more accurate inference,
also based on the likelihood function. In Section~\ref{sec4}, I~describe some
extensions of the likelihood function that have been proposed for
models with complex dependence structure, with particular emphasis on
composite likelihood.

%s2 ###
\section{Inference based on the likelihood function}\label{sec2}

Suppose we have a probability model for an observable random vector $Y
= (Y_1, \ldots, Y_n)$ of the form $f(y;\theta)$, where $\theta$ is a
vector of unknown parameters in the model, and $f(y; \theta)$ is a
density function with respect to a dominating measure, usually Lebesgue
measure or counting measure, depending on whether our observations are
discrete or continuous.
%discrete data, on the grounds that observables cannot be measured with
%infinite precision, but the use of continuous models as a mathematical
%approximation seems to lead to a much simpler theory.}
Typical models used in applications assume that $\theta$ could
potentially be any value in a set~$\Theta$; sometimes $\Theta$ is
infinite-dimensional, but more usually $\Theta\subset\mathbb{R}^d$.
The inverse problem mentioned in Section~\ref{sec1} is to construct inference
about the value or values of $\theta\in\Theta$ that could plausibly
have generated an observed value $y = y^0$. This is a considerable
abstraction from realistic applied settings; in most scientific work
such a problem will not be isolated from a series of investigations,
but we can address at least some of the main issues in this setting.

The likelihood function is simply
%
%e2.1 ###
\begin{equation}
\label{likelihood} L(\theta; y) \propto f(y;\theta);
\end{equation}
i.e., there is an equivalence class of likelihood functions $L(\theta
;y) = c(y)f(y;\theta)$, and only relative ratios $L(\theta_2;y)/L(\theta_1;y)$ are uniquely determined. From a mathematical
point of view, (\ref{likelihood}) is a trivial re-expression of the
model $f(y;\theta)$; the re-ordering of the arguments is simply to
emphasize in the notation that we are more interested in the $\theta
$-section for fixed $y$ than in the $y$-section for fixed $\theta$.
Used directly with a given observation $y^0$, $L(\theta;y^0)$ provides
a ranking of relative plausibility of various values of $\theta$, in
light of the observed data.

A form of direct inference can be obtained by plotting the likelihood
function, if the parameter space is one- or two-dimensional, and
several writers, including Fisher, have suggested declaring values of
$\theta$ in ranges determined by likelihood ratios as plausible, or
implausible; for example, \citet{fisher56} suggested that values of
$\theta$ for which $L(\hat\theta;y)/L(\theta;y) > 15$, be declared
`implausible', where $\hat\theta=\hat\theta(y)$ is the maximum
likelihood estimate of $\theta$, i.e., the value for which the
likelihood function is maximized, over $\theta$, for a given $y$.

In general study of statistical theory and methods we are usually
interested in properties of our statistical methods, in repeated
sampling from the model $f(y;\theta_0)$, where $\theta_0$ is the
notional `true' value of $\theta$ that generated the data. This
requires considering the distribution of $L(\theta;Y)$, or relative
ratios such as $L\{\hat\theta(Y);Y\}/L\{\theta(Y);Y\}$. To this end,
some standard summary functions of $L(\theta;Y)$ are defined. Writing
$\ell(\theta;Y) = \log L(\theta; Y)$ we define the \emph{score
function} $u(\theta;Y) = \partial\ell(\theta;Y)/\partial\theta$,
and the observed and expected Fisher information functions:
%
%e2.2 ###
\begin{equation}
\label{fisherinfo} j(\theta;Y) = - \frac{\partial^2\ell(\theta;Y)}{\partial\theta\,
\partial\theta^T},\qquad i(\theta) = E\biggl
\{- \frac{\partial^2\ell(\theta
;Y)}{\partial\theta\,\partial\theta^T}\biggr\}.
\end{equation}

If the components of $Y$ are independent, then $\ell(\theta;Y)$ is a
sum of independent random variables, as is $u(\theta;Y)$, and under
some conditions on the model the central limit theorem for $u(\cdot
;Y)$ leads to the following asymptotic results, as $n\rightarrow\infty$:
%
%e2.5 ###
%e2.4 ###
%e2.3 ###
\begin{eqnarray}
\label{score} s(\theta) &=& j^{-1/2}(\hat\theta)u(\theta) \goesd N(0,
I),
\\
\label{mle}q(\theta) &=& j^{1/2}(\hat\theta) (\hat\theta-\theta) \goesd
N(0, I),
\\
\label{lrt}w(\theta) &=& 2\bigl\{\ell(\hat\theta)-\ell(\theta)\bigr\} \goesd
\chi^2_d,
\end{eqnarray}
where we suppress the dependence of each derived quantity on $Y$ (and
on $n$) for notational convenience. These results hold under the model
$f(y;\theta)$; a more precise statement would use the true value
$\theta_0$ in $u(\theta)$, $(\hat\theta-\theta)$, and $\ell
(\theta)$ above, and the model $f(y;\theta_0)$. However, the
quantities $s(\theta)$, $q(\theta)$ and $w(\theta)$, considered as
functions of both $\theta$ and $Y$, are approximate \emph{pivotal
quantities}, i.e., they have a known distribution, at least
approximately. For $\theta\in\mathbb{R}$ we could plot, for example,
$\Phi\{q(\theta)\}$ as a function of $\theta$, where $\Phi(\cdot)$
is the standard normal distribution function, and obtain approximate
$p$-values for testing any value of $\theta\in\mathbb{R}$ for fixed
$y$. The approach to inference based on these pivotal quantities avoids
the somewhat artificial distinction between point estimation and
hypothesis testing. When $\theta\in\mathbb{R}$, an approximately
standard normal pivotal quantity can be obtained from (\ref{lrt}) as
%
%e2.6 ###
\begin{equation}
\label{root} r(\theta) = \operatorname{sign}(\hat\theta-\theta)\bigl[2\bigl\{
\ell(\hat\theta )-\ell(\theta)\bigr\}\bigr]^{1/2} \goesd N(0, 1).
\end{equation}

The likelihood function is also the starting point for Bayesian
inference; if we model the unknown parameter as a random quantity with
a postulated prior probability density function $\pi(\theta)$, then
inference given an observed value $Y=y$ is based on the posterior
distribution, with density
%
%e2.7 ###
\begin{equation}
\label{bayes} \pi(\theta\mid y) =\frac{ \exp\{\ell(\theta;y)\}\pi(\theta
)}{\int\exp\{\ell(\phi;y)\}\pi(\phi)\,d\phi}.
\end{equation}
Bayesian inference is conceptually straightforward, given a prior
density, and computational methods for estimating the integral in the
denominator of (\ref{bayes}), and associated integrals for marginal
densities of components, or low-dimensional functions of $\theta$,
have enabled the application of Bayesian inference in models of
considerable complexity. Two very useful methods include Laplace
approximation of the relevant integrals, and Markov chain Monte Carlo
simulation from the posterior. Difficulties with Bayesian inference
include the specification of a prior density, and the meaning of
probabilities for parameters of a mathematical model.

One way to assess the influence of the prior is to evaluate the
properties of the resulting inference under the sampling model, and
under regularity conditions similar to those needed to obtain (\ref
{score}), (\ref{mle}) and (\ref{lrt}), a normal approximation to the
posterior density can be derived:
%
%e2.8 ###
\begin{equation}
\label{bayes2} \pi(\theta\mid y) \dotsim N\bigl\{\hat\theta, j^{-1}(
\hat\theta)\bigr\},
\end{equation}
implying that inferences based on the posterior are asymptotically
equivalent to those based on $q$. This simple result underlines the
fact that Bayesian inference will in large samples give approximately
correct inference under the model, and also that to distinguish between
Bayesian and non-Bayesian approaches we need to consider the next order
of approximation.

If $\theta\in\mathbb{R}^d$, then (\ref{score})--(\ref{lrt}) can
be used to construct confidence regions, or to test simple hypotheses
of the form $\theta= \theta_0$, but in many settings $\theta$ can
usefully be separated into a parameter of interest $\psi$, and a
nuisance parameter $\lambda$, and analogous versions of the above
limiting results in this context are
%
%e2.11 ###
%e2.10 ###
%e2.9 ###
\begin{eqnarray}
\label{score2}s(\psi) &=& j_{\mathrm{p}}^{-1/2}(\hat\psi)
\ell_{\mathrm{p}}'(\psi) \goesd N(0, I),
\\
\label{mle2}q(\psi) &=& j_{\mathrm{p}}^{1/2}(\hat\psi) (\hat\psi-
\psi) \goesd N(0, I),
\\
\label{lrt2}w(\psi) &=& 2\bigl\{\ell_{\mathrm{p}}(\hat\psi) -
\ell_{\mathrm{p}}(\psi)\bigr\} \goesd\chi^2_{d_1},
\end{eqnarray}
where $\ell_{\mathrm{p}}(\psi) = \ell(\psi, \hat\lambda_\psi)$ is
the profile log-likelihood function, $\hat\lambda_\psi$ is the
constrained maximum likelihood estimate of the nuisance parameter
$\lambda$ when $\psi$ is fixed, $d_1$ is the dimension of $\psi$,
and $j_{\mathrm{p}}(\psi) = -\partial^2\ell_{\mathrm{p}}(\psi)/\partial
\psi\,\partial\psi^T$ is the Fisher information function based on the
profile log-likelihood function.

The third result (\ref{lrt2}) can be used for model assessment among
nested models; for example, the exponential distribution is nested
within both the Gamma and Weibull models, and a test based on $w$ of,
say, a gamma model with unconstrained shape parameter, and one with the
shape parameter set equal to 1, is a test of fit of the exponential
model to the data; the rate parameter is the nuisance parameter
$\lambda$. The use of the log-likelihood ratio to compare two
non-nested models, for example a log-normal model to a gamma model,
requires a different asymptotic theory (\citeauthor{ch}, \citeyear{ch}, Ch. 8).

A related approach to model selection is based on the Akaike
information criterion,
\[
\mathit{AIC} = -2\ell(\hat\theta) + 2d,
\]
where $d$ is the dimension of $\theta$. Just as only differences in
log-likelihoods are relevant, so are differences in $\mathit{AIC}$: for a
sequence of model fits the one with the smallest value of $\mathit{AIC}$ is
preferred. The $\mathit{AIC}$ criterion was developed in the context of
prediction in time series, but can be motivated as an estimate of the
Kullback-Leibler divergence between a fitted model and a notional
`true' model. The statistical properties of $\mathit{AIC}$ as a model selection
criterion depend on the context; for example for choosing among a
sequence of regression models of the same form, model selection using
$\mathit{AIC}$ is not consistent (\citeauthor{dav}, \citeyear{dav}, Ch.~4.7).\vadjust{\goodbreak}
Several related versions of model selection criterion have been suggested, including
modifications to $\mathit{AIC}$, and a version motivated by Bayesian arguments,
\[
\mathit{BIC} = -2\ell(\hat\theta) + d\log(n),
\]
where $n$ is the sample size for the model with $d$ parameters.

%s3 ###
\section{More accurate inference}\label{sec3}

The approximate inference suggested by the approximate pivotal
quantities (\ref{score2}), (\ref{mle2}) and (\ref{lrt2}) is obtained
by treating the profile log-likelihood function as if it were a genuine
log-likelihood function, i.e. as if the true value of $\lambda$ were
$\hat\lambda_{\psi}$. This can be misleading, because it does not
account for the fact that the nuisance parameter has been estimated.
One familiar example is inference for the variance in a normal theory
linear regression model; the maximum likelihood estimate is
\[
\hat\sigma^2 = (y-X\hat\beta)^T(y-X\hat\beta)/n,
\]
which has expectation $(n-k)\sigma^2/n$, where $k$ is the dimension of
$\beta$. Although this estimator is consistent as $n\rightarrow\infty
$ with $k$ fixed, it can be a poor estimate for finite samples,
especially if $k$ is large relative to $n$, and the divisor $n-k$ is
used in practice. One way to motivate this is to note that $n\hat
\sigma^2/(n-k)$ is unbiased for $\sigma^2$; an argument that
generalizes more readily is to note that the likelihood function
$L(\beta,\sigma^2;\hat\beta,\hat\sigma^2)$ can be expressed as
\[
L_1\bigl(\mu, \sigma^2; \hat\beta\bigr) L_2
\bigl(\sigma^2;\hat\sigma^2\bigr),
\]
where $L_1$ is proportional to the density of $\hat\beta$ and $L_2$
is the marginal density of $\hat\sigma^2$ or equivalently $(y-X\hat
\beta)^T(y-X\hat\beta)$. The unbiased estimate of $\sigma^2$
maximizes the second component $L_2$, which is known as the restricted
likelihood, and estimators based on it often called ``REML'' estimators.

Higher order asymptotic theory for likelihood inference has proved to
be very useful for generalizing these ideas, by refining the profile
log-likelihood to take better account of the nuisance parameter, and
has also provided more accurate distribution approximations to pivotal
quantities. Perhaps most importantly, for statistical theory, higher
order asymptotic theory helps to
clarify the role of the likelihood function and the prior in the
calibration of Bayesian inference.
These three goals have turned out to be very intertwined.%\cite{r1}

To illustrate some aspects of this, consider the marginal posterior
density for $\psi$, where $\theta= (\psi, \lambda)$:
%
%e3.1 ###
\begin{equation}
\pi_m(\psi\mid y) = \frac{\int\exp\{\ell(\psi,\lambda;y)\}\pi
(\psi,\lambda)\,d\lambda}{\int\exp\{\ell(\psi,\lambda;y)\}\pi
(\psi,\lambda)\,d\lambda \,d\psi}.
\end{equation}
Laplace approximation to the numerator and denominator integrals leads to
%
%e3.2 ###
\begin{eqnarray}
\label{laplace2}\pi_m(\psi\mid y) &\doteq& \frac{(2\pi)^{(d-d_1)/2} \exp\{\ell
(\psi,\hat\lambda_\psi)\}|j_{\lambda\lambda}(\psi,\hat\lambda_\psi)|^{-1/2}\pi(\psi,\hat\lambda_\psi)} {
(2\pi)^{d/2} \exp\{\ell(\hat\psi,\hat\lambda)\}|j(\hat\psi
,\hat\lambda)|^{-1/2}\pi(\hat\psi,\hat\lambda)}
\\
&=& \frac{1}{(2\pi)^{d_1/2}}\exp\bigl\{\ell_{\mathrm{p}}(\psi)-\ell_{\mathrm{p}}(
\hat\psi)\bigr\}|j_{\mathrm{p}}(\hat\psi)|^{1/2} \biggl\{
\frac
{|j_{\lambda\lambda}(\psi,\hat\lambda_\psi)|}{|j_{\lambda\lambda
}(\hat\psi,\hat\lambda)|} \biggr\}^{-1/2}\frac{\pi(\psi,\hat
\lambda_\psi)}{\pi(\hat\psi,\hat\lambda)},
\nonumber
\\
&=& \frac{1}{(2\pi)^{d_1/2}}\exp\bigl\{\ell_{\mathrm{a}}(\psi) - \ell_{\mathrm{a}}(
\hat\psi)\bigr\} \frac{\pi(\psi,\hat\lambda_\psi)}{\pi
(\hat\psi,\hat\lambda)},
\nonumber
\end{eqnarray}
where $j_{\lambda\lambda}(\theta)$ is the block of the observed
Fisher information function corresponding to the nuisance parameter
$\lambda$, $|j(\hat\theta)|$ has been computed using the partitioned
form to give the second expression in (\ref{laplace2}),
and in the third expression
\[
\ell_{\mathrm{a}}(\psi) = \ell_{\mathrm{p}}(\psi) - (1/2)
\log\bigl|j_{\lambda
\lambda}(\psi,\hat\lambda_\psi)\bigr|.
\]
When renormalized to integrate to one, this Laplace approximation has
relative error $O(n^{-3/2})$ in independent sampling from a model that
satisfies various regularity conditions similar to those needed to show
the asymptotic normality of the posterior \cite{tk}.

These expressions show that an adjustment for estimation of the
nuisance parameter is captured in $\log|j_{\lambda\lambda}(\cdot
)|$, and this adjustment can be included in the profile log-likelihood
function, as in the third expression in (\ref{laplace2}), or tacked
onto it, as in the second expression. The effect of the prior is
isolated from this nuisance parameter adjustment effect, so, for
example, if $\hat\lambda_\psi= \hat\lambda$, and the priors for
$\psi$ and $\lambda$ are independent, then the form of the prior for
$\lambda$ given $\psi$ does not affect the approximation.

The adjusted profile log-likelihood function $\ell_{\mathrm{a}}(\psi)$
is the simplest of a number of modified profile log-likelihood
functions suggested in the literature for improved frequentist
inference in the presence of nuisance parameters, and was suggested for
general use in \citet{cr}, after reparametrizing the model to make
$\psi$ and $\lambda$ orthogonal with respect to expected Fisher
information, i.e., ${\mathrm{E}}\{-\partial^2\ell(\psi,\lambda
)/\partial\psi\,\partial\lambda\} = 0$. This reparameterization makes
it at least more plausible that $\psi$ and $\lambda$ could be
modelled as \emph{a priori} independent, and also ensures that $\hat
\lambda_{\psi} - \hat\lambda= O_p(1/n)$, rather than the usual
$O_p(1/\surd n)$.

A number of related, but more precise, adjustments to the profile
log-likelihood function have been developed from asymptotic expansions
for frequentist inference, and take the form
%
%e3.3 ###
\begin{equation}
\label{mpl} \ell_{\mathrm{M}}(\psi) = \ell_{\mathrm{p}}(\psi) + (1/2)
\log|j_{\lambda
\lambda}(\psi, \hat\lambda_{\psi})| + B(\psi),
\end{equation}
where $B(\psi) = O_p(1)$;
see, for example, \citeauthor{dmsy} (\citeyear{dmsy}) and \citet{ps}. The change from $-1/2$
to $+1/2$ is related to the orthogonality conditions; in (\ref{mpl})
orthogonality of parameters is not needed, as the expression is
parameterization invariant.

Inferential statements based on approximations from
(\ref{score2})--(\ref{lrt2}), with $\ell_{\mathrm{a}}(\psi)$ or $\ell_{\mathrm{M}}(\psi
)$ substituting for the profile log-likelihood function, are still
valid and are more accurate in finite samples, as they adjust for
errors due to estimation of $\lambda$. They are still first-order
approximations, although often quite good ones.

One motivation for these modified profile log-likelihood functions, and
inference based on them, is that they approximate marginal or
conditional likelihoods, when these exist. For example, if the model is
such that
\[
f(y;\psi, \lambda) \propto g_1(t_1;\psi)
g_2(t_2\mid t_1;\lambda),
\]
then inference for $\psi$ can be based on the marginal likelihood for
$\psi$ based on $t_1$, and the theory outlined above applies directly.
This factorization is fairly special; more common is a factorization of
the form $g_1(t_1;\psi) g_2(t_2\mid t_1;\lambda, \psi)$: in that
case to base our inference on the likelihood for $\psi$ from $t_1$
would require further checking that little information is lost in
ignoring the second term. Arguments like these, applied to special
classes of model families, were used to derive the modified profile
log-likelihood inference outlined above.

A related development is the improvement of the distributional
approximation to the approximate pivotal quantity (\ref{root}). The
Laplace approximation (\ref{laplace2}) can be used to obtain the
Bayesian pivotal, for scalar $\psi$,
%
%e3.4 ###
\begin{equation}
\label{rstarb} r^*_B(\psi) = r(\psi) + \frac{1}{r(\psi)}\log\biggl
\{\frac{q_B(\psi
)}{r(\psi)}\biggr\} \dotsim N(0,1),
\end{equation}
where
%
%e3.6 ###
%e3.5 ###
\begin{eqnarray}
\label{root2}r(\psi) &=& \operatorname{sign}(\hat\psi- \psi)\bigl[2\bigl\{
\ell_{\mathrm{p}}(\hat \psi) - \ell_{\mathrm{p}}(\psi)\bigr\}
\bigr]^{1/2},
\\
\label{qb}q_B(\psi) &=& -\ell_{\mathrm{p}}'(\psi)
j_{\mathrm{p}}^{-1/2}(\hat\psi) \biggl\{\frac{|j_{\lambda\lambda}(\psi,\hat\lambda_\psi
)|}{|j_{\lambda\lambda}(\hat\psi,\hat\lambda)|} \biggr
\}^{1/2}\frac{\pi(\hat\psi,\hat\lambda)}{\pi(\psi,\hat\lambda_\psi)}
\end{eqnarray}
and the approximation in (\ref{rstarb}) is to the posterior
distribution of $r^*$, given $y$, and is accurate to $O(n^{-3/2})$.

There is a frequentist version of this pivotal that has the same form:
%
%e3.7 ###
\begin{equation}
\label{rstarf} r^*_F(\psi) = r(\psi) + \frac{1}{r(\psi)}\log\biggl
\{\frac{q_F(\psi
)}{r(\psi)}\biggr\} \dotsim N(0,1),
\end{equation}
where
$r(\psi)$ is given by (\ref{root2}), but the expression for $q_F(\psi
)$ requires additional notation, and indeed an additional likelihood
component. In the special case of no nuisance parameters
%
%e3.9 ###
%e3.8 ###
\begin{eqnarray}
\label{qf1}q_F(\theta) &=& \bigl\{\ell_{;\hat\theta}(\hat\theta;
\hat\theta, a)-\ell_{;\hat\theta}(\theta;\hat\theta,a)\bigr\} j^{-1/2}(
\hat \theta;\hat\theta,a)
\\
\label{qf2}&=&\bigl\{\varphi(\hat\theta)-\varphi(\theta)\bigr\}
\varphi_\theta^{-1}(\hat\theta)j^{1/2}(\hat\theta).
\end{eqnarray}
In (\ref{qf1}), we have assumed that there is a one-to-one
transformation from $y$ to $(\hat\theta, a)$, and that we can write
the log-likelihood function in terms of $\theta, \hat\theta,a$ and
then differentiate it with respect to $\hat\theta$, for fixed $a$.
Expression (\ref{qf2}) is equivalent, but expresses this sample space
differentiation through a data-dependent reparameterization $\varphi
(\theta) = \varphi(\theta;y) = \partial\ell(\theta;y)/\partial
V(y)$, where the derivative with respect to $V(y)$ is a directional
derivative to be determined.

The details are somewhat cumbersome, and even more so for the case of
nuisance parameters, but the resulting $r^*_F$ approximate pivotal
quantity is readily calculated in a wide range of models for
independent observations $y_1, \ldots, y_n$. Detailed accounts are
given in \citet{bnc}, \citet{psbook}, \citet{sevbook}, \citet{frw}
and \citeauthor{bdr} (\citeyear{bdr}, Ch. 8.6); the last emphasizes implementation in a
number of practical settings, including generalized linear models,
nonlinear regression with normal errors, linear regression with
non-normal errors, and a number of more specialized models.

From a theoretical point of view, an important distinction between
$r^*_B$ and $r^*_F$ is that the latter requires differentiation of the
log-likelihood function on the sample space, whereas the former depends
only on the observed log-likelihood function, along with the prior. The
similarity of the two expressions suggests that it might be possible to
develop prior densities for which the posterior probability bounds are
guaranteed to be valid under the model, at least to a higher order of
approximation than implied by (\ref{bayes2}), and there is a long line
of research on the development of these so-called ``matching priors'';
see, for example, \citet{dmuk}.

%s4 ###
\section{Extending the likelihood function}\label{sec4}

%s4.1 ###
\subsection{Introduction}

While the asymptotic results of the last section provide very accurate
inferences, they are not as straightforward to apply as the first order
results, especially in models with complex dependence. They do shed
light on many aspects of theory, including the precise points of
difference, asymptotically, between Bayesian and nonBayesian inference.
And the techniques used to derive them, saddlepoint and Laplace
approximations in the main, have found application in complex models in
certain settings, such as the integrated nested Laplace approximation
of \citet{rmc}.

A glance at any number of papers motivated by specific applications,
though, will confirm that likelihood summaries, and in particular
computation of the maximum likelihood estimator, are often the
inferential goal, even as the models become increasingly high-dimensional.

This is perhaps a natural consequence of the emphasis on developing
probability models that could plausibly generate, or at least describe,
the observed responses, as the likelihood function is directly obtained
from the probability model. But more than this, inference based on the
likelihood function provides a standard set of tools, whose properties
are generally well-known, and avoids the construction of \emph{ad hoc}
inferential techniques for each new application. For example,
%in the field of computational neuroscience
\citet{feng} write ``The likelihood framework is an efficient way to
extract information from a neural spike train\ldots We believe that
greater use of the likelihood based approaches and goodness-of-fit
measures can help improve the quality of neuroscience data analysis''.

A number of inference functions based on the likelihood function, or
meant to have some of the key properties of the likelihood function,
have been developed in the context of particular applications or
particular model families. In some cases the goal is to find
`reasonably reliable' estimates of a parameter, along with an estimated
standard error; in other cases the goal is to use approximate pivotal
quantities like those outlined in Section~\ref{sec2} in settings where the
likelihood is difficult to compute. The goal of obtaining reliable
likelihood-based inference in the presence of nuisance parameters was
addressed in Section~\ref{sec3}. In some settings, families of parametric models
are too restrictive, and the aim is to obtain likelihood-type results
for inference in semi-parametric and non-parametric settings.

%s4.2 ###
\subsection{Generalized linear mixed models}

In many applications with longitudinal, clustered, or spatial data, the
starting point is a generalized linear model with a linear predictor of
the form $X\beta+ Zu$, where $X$ and $Z$ are $n\times k$ and $n\times
q$, respectively, matrices of predictors, and $u$ is a $q$-vector of
random effects. The marginal distribution of the responses requires
integrating over the distribution of the random effects $u$, and this
is often computationally infeasible. Many approximations have been
suggested: one approach is to approximate the integral by Laplace's
method \cite{breslow}, leading to what is commonly called penalized
quasi-likelihood, although this is different from the penalized
versions of composite likelihood discussed below. The term
quasi-likelihood in the context of generalized linear models refers to
the specification of the model through the mean function and variance
function only, without specifying a full joint density for the
observations. This was first suggested by \citet{wedderburn}, and
extended to longitudinal data in \citet{liang} and later work, leading
to the methodology of generalized estimating equations, or GEE. \citet
{renard} compared penalized quasi-likelihood to pairwise likelihood,
discussed in Section~\ref{sec4.3}, in simulations of multivariate probit models
for binary data with random effects. In general penalized
quasi-likelihood led to estimates with larger bias and variance than
pairwise likelihood.

A different approach to generalized linear mixed models has been
developed by Lee and Nelder; see, for example, \citet{leenelder} and
\citet{nelderlee}, under the name of $h$-likelihood. This addresses
some of the failings of the penalized quasi-likelihood method by
modelling the mean parameters and dispersion parameters separately. The
$h$-likelihood for the dispersion parameters is motivated by REML-type
arguments not unrelated to the higher order asymptotic theory outlined
in the previous section. There are also connections to work on
prediction using likelihood methods \cite{bjornstadt}. Likelihood
approaches to prediction have proved to be somewhat elusive, at least
in part because the `parameter' to be predicted is a random variable,
although Bayesian approaches are straightforward as no distinction is
made between parameters and random variables.

%s4.3 ###
\subsection{Composite likelihood}\label{sec4.3}

Composite likelihood is one approach to combining the advantages of
likelihood with computational feasibility; more precisely it is a
collection of approaches. The general principle is to simplify complex\vadjust{\goodbreak}
dependence relationships by computing marginal or conditional
distributions of some subsets of the responses, and multiplying these
together to form an inference function.

As an \emph{ad hoc} solution it has emerged in several versions and in
several contexts in the statistical literature; an important example is
the pseudo-likelihood for spatial processes proposed in \citeauthor{besag1}
(\citeyear{besag1,besag2}). In studies of large networks, computational complexity can be
reduced by ignoring links between distant nodes, effectively treating
sub-networks as independent. In Gaussian process models with
high-dimensional covariance matrices, assuming sparsity in the
covariance matrix is effectively assuming subsets of variables are
independent. The term composite likelihood was proposed in \citet
{lindsay88}, where the theoretical properties of composite likelihood
estimation were studied in some generality.

We suppose a vector response of length $q$ is modelled by $f(y;\theta
), \theta\in\mathbb{R}^d$. Given a set of events $\mathcal{A}_k, k
= 1, \ldots, K$, the composite likelihood function is defined as
%
%e4.1 ###
\begin{equation}
\label{cl} \mathit{CL}(\theta;y) = \prod_{k=1}^K
f(y \in\mathcal{A}_k; \theta),
\end{equation}
and the composite log-likelihood function is
%
%e4.2 ###
\begin{equation}
\label{logcl} \mathit{c\ell}(\theta;y) = \sum_k \log f( y
\in\mathcal{A}_k;\theta).
\end{equation}
Because each component in the sum is the log of a density function, the
resulting score function $\partial \mathit{c\ell}(\theta;y)/\partial\theta$
has expected value $0$, so has at least one of the properties of a
genuine log-likelihood function.

Relatively simple and widely used examples of composite likelihoods
include independence composite likelihood,
\[
\mathit{c\ell}_{\mathrm{ind}}(\theta;y) = \sum_{r=1}^q
\log f_1(y_r;\theta),
\]
pairwise composite likelihood
\[
\mathit{c\ell}_{\mathrm{pair}}(\theta;y) = \sum_{r=1}^q
\sum_{s>r} \log f_2(y_r,
y_{s};\theta),
\]
and pairwise conditional composite likelihood
%
%e4.3 ###
\begin{equation}
\label{clcond} \mathit{c\ell}_{\mathrm{cond}}(\theta;y) = \sum
_{r=1}^q \log f(y_r\mid
y_{(-r)};\theta),
\end{equation}
where $f_1(y_r;\theta)$ and $f_2(y_r, y_s;\theta)$ are the marginal
densities for a single component and a pair of components of the vector
observation, and the density in (\ref{clcond}) is the conditional
density of one component, given the remainder.

Many similar types of composite likelihood can be constructed,
appropriate to time series, or spatial data, or repeated measures, and
so on, and the definition\vadjust{\goodbreak} is usually further extended by allowing each
component event to have an associated weight $w_k$. Indeed one of the
difficulties of studying the theory of composite likelihood is the
generality of the definition.

Inference based on composite likelihood is constructed from analogues
to the asymptotic results for genuine likelihood functions. Assuming we
have a sample $\underline y = (y^{(1)}, \ldots, y^{(n)})$ of
independent observations of $y$, the composite score function,
%
%e4.4 ###
\begin{equation}
\label{clscore} u_{\mathit{CL}}(\theta;\underline y) = \sum
_{i=1}^n \sum_k
\partial\log f\bigl(y^{(i)} \in\mathcal{A}_k;\theta\bigr)/
\partial\theta,
\end{equation}
is used as an estimating function to obtain the maximum composite
likelihood estimator $\hat\theta_{\mathit{CL}}$, and under regularity
conditions on the full model, with $n \rightarrow\infty$ and fixed
$K$, we have, for example,
%
%e4.5 ###
\begin{equation}
\label{cmle} (\hat\theta_{\mathit{CL}} - \theta)^TG(\hat
\theta_{\mathit{CL}}) (\hat\theta_{\mathit{CL}}-\theta) \goesd
\chi^2_d,
\end{equation}
where
%
%e4.6 ###
\begin{equation}
\label{godambe} G(\theta) = H(\theta)J^{-1}(\theta)H(\theta)
\end{equation}
is the $d\times d$ Godambe information matrix, and
\[
J(\theta) = \operatorname{var}\bigl\{u_{\mathit{CL}}(\theta;Y)\bigr\},\qquad H(
\theta) = \mathrm{E}\bigl\{ -(\partial/\partial\theta) u_{\mathit{CL}}(\theta;Y)
\bigr\},
\]
are the variability and sensitivity matrix associated with $u_{\mathit{CL}}$.

The analogue of (\ref{lrt}) is
%
%e4.7 ###
\begin{equation}
\label{clrt} 2\bigl\{\mathit{c\ell}(\hat\theta_{\mathit{CL}})-\mathit{c\ell}(\theta)\bigr\} \goesd
\sum_{i=1}^d \lambda_i
\chi^2_{1i},
\end{equation}
where $\lambda_i$ are the eigenvalues of $J^{-1}(\theta)H(\theta)$.

Neither of these results is quite as convenient as the full likelihood
versions, and in particular contexts it may be difficult to estimate
$J(\theta)$ accurately, but there are a number of practical settings
where these results are much more easily implemented than the full
likelihood results, and the efficiency of the methods can be quite good.

A number of applied contexts are surveyed in \citet{vrf}. As just one
example, developed subsequently, \citet{dpr} investigate pairwise
composite likelihood for max-stable processes, developed to model
extreme values recorded at a number $D$ of spatially correlated sites.
Although the form of the $D$-dimensional density is known, it is not
computable for $D>3$, although expressions are available for the joint
density at each pair of sites. Composite likelihood seems to be
particularly important for various types of spatial models, and many
variations of it have been suggested for these settings.

In some applications, particularly for time series, but also for
space-time data, a~sample of independent observations is not available,
and the relevant asymptotic theory is for $q \rightarrow\infty$,
where $q$ is the dimension of the single response. The asymptotic
results outlined above will require some conditions on the decay of the
dependence among components as the `distance' between them increases.
Asymptotic theory for pairwise likelihood is investigated in \citet
{davisyau} for linear time series, and in \citet{davisklupp} for
max-stable processes in space and time.

Composite likelihood can also be used for model selection, with an
expression analogous to $\mathit{AIC}$, and for Bayesian inference, after
adjustment to accommodate result (\ref{clrt}). \emph{Statistica
Sinica} \textbf{21}, \#1 is a special issue devoted to composite
likelihood, and more recent research is summarized in the report on a
workshop at the Banff International Research Station \cite{joe}.

%s4.4 ###
\subsection{Semi-parametric likelihood}

In some applications, a flexible class of models can be constructed in
which the nuisance `parameter' is an unknown function. The most
widely-known example is the proportional hazards model of \citet
{cox72} for censored survival data; but semi-parametric regression
models are also widely used, where the particular covariates of
interest are modelled with a low-dimensional regression parameter, and
other features expected to influence the response are modelled as
`smooth' functions. \citet{cox72} developed inference based on a
partial likelihood, which ignored the aspects of the likelihood bearing
on the timing of failure events, and subsequent theory based on
asymptotics for counting processes established the validity of this
approach. In fact, \citet{cox72}'s partial likelihood can be viewed as
an example of composite likelihood as described above, although the
theory for general semi-parametric models seems more natural.

\citet{mvv} showed that partial likelihood can be viewed as a profile
likelihood, maximized over the nuisance function, and discussed a class
of semi-parametric models for which the profile likelihood continues to
have the same asymptotic properties as the usual parametric profile
likelihood; the contributions to the discussion of their results
provide further insight and references to the extensive literature on
semi- and non-parametric likelihoods. There is, however, no guarantee
that asymptotic theory will lead to accurate approximation for finite
samples; it would presumably have at least the same drawbacks as
profile likelihood in the parametric setting. Improvements via
modifications to the profile likelihood, as described above in the
parametric case, do not seem to be available in these more general settings.

Some semi-parametric models are in effect converted to high-dimensional
parametric models through the use of linear combinations of basis
functions; thus the linear predictor associated with a component $y_i$
might be $\beta_0 + \beta_1 x_i + \sum_{j=1}^J \gamma_jB(z_i)$, or
$\beta_0 + \beta_1 x_i + \sum_{j=1}^J \gamma_{1j}B_j(z_{1i}) +
\cdots+ \sum_{j=1}^J \gamma_{kj}B_j(z_{ki})$. The log-likelihood
function for models such as these is often regularized, so that
$\ell(\beta,\gamma)$ is replaced by $\ell(\beta,\gamma) + \lambda
p(\gamma)$, where $p(\cdot)$ is a penalty function such as $\Sigma
\gamma_{kj}^2$ or $\Sigma|\gamma_{kj}|$, and $\lambda$ a tuning parameter.
Many of these extensions, and the asymptotic theory associated with
them, are discussed in \citeauthor{vdv} (\citeyear{vdv}, Ch. 25). Penalized likelihood using
squared error is reviewed in \citet{green}; the $L_1$ penalty has been
suggested as a means of combining likelihood inference with variable
selection; see, for example, \citet{fan}.

Penalized composite likelihoods have been proposed for applications in
spatial analysis (\citeauthor{divino}, \citeyear{divino};
\citeauthor{apanasovich}, \citeyear{apanasovich};
\citeauthor{xue}, \citeyear{xue}), Gaussian graphical
models \cite{gaomassam}, and clustered longitudinal data \cite{gaosong}.

The difference between semi-parametric likelihoods and nonparametric
likelihoods is somewhat blurred; both have an effectively
infinite-dimensional parameter space, and as discussed in \citet{mvv}
and the discussion, conditions on the model to ensure that
likelihood-type asymptotics still hold can be quite technical.

Empirical likelihood is a rather different approach to non-parametric
models first proposed by \citet{owen}; a recent discussion is
\citet{hmvk}. Empirical likelihood assumes the existence of a
finite-dimensional parameter of
interest, defined as a functional of the distribution function for the
data, and constructs a profile likelihood by maximizing the joint
probability of the data, under the constraint that this parameter is
fixed. This construction is particularly natural in survey sampling,
where the parameter is often a property of the population (\citeauthor{cs}, \citeyear{cs};
\citeauthor{wr}, \citeyear{wr}). Distribution theory for empirical likelihood more closely follows
that for usual parametric likelihoods.

%s4.5 ###
\subsection{Simulation methods}

Simulation of the posterior density by Markov chain Monte Carlo methods
is widely used for Bayesian inference, and there is an enormous
literature on various methods and their properties. Some of these
methods can be adapted for use when the likelihood function itself
cannot be computed, but it is possible to simulate observations from
the stochastic model; many examples arise in statistical genetics.
Simulation methods for maximum likelihood estimation in genetics was
proposed in \citet{gt}; more recently sequential Monte Carlo methods
(see, for example, \citet{sisson}) and ABC (approximate Bayesian
computation) methods (\citeauthor{fp}, \citeyear{fp}; \citeauthor{mprr}, \citeyear{mprr}) are being investigated as
computational tools.

%s5 ###
\section{Conclusion}

A reviewer of an earlier draft suggested that a great many
applications, especially involving very large and/or complex datasets,
take more algorithmic approaches, often using techniques designed to
develop sparse solutions, such as wavelet or thresholding techniques,
and that likelihood methods may not be relevant for these
application areas.

Certainly a likelihood-based approach depends on a statistical
model for the data, and for many applications under the general rubric
of machine learning these may not be as important as developing fast
and reliable approaches to prediction; recommender systems are one such example.

There are however many applications of `big data'
methods where statistical models do provide some structure, and in
these settings, as in the more classical application areas, likelihood
methods provide a unifying basis for inference.

\section*{Acknowledgements}
This research was partially supported by the Natural Sciences and
Engineering Research Council. Thanks are due to two reviewers for
helpful comments on an earlier version.

%suskaldyti doi

% imsref loaded by audrone.aklyte, 2012-12-18 10:57:48
%
% imsref loaded by audrone.aklyte, 2012-12-18 11:06:32
% imsref loaded by audrone.aklyte, 2012-12-18 13:18:40
% imsref loaded by audrone.aklyte, 2012-12-18 16:12:11

\end{document}